\documentclass[a4paper,12pt]{article}
\usepackage{amsthm,amsmath,amssymb,txfonts,graphics,,mathrsfs}
\usepackage[latin1]{inputenc}

\usepackage{pb-diagram}
\usepackage{hyperref}

\newtheorem{lemma}{Lemma}
\newtheorem{teo}{Theorem}
\newtheorem{oss}{Remark}
\newtheorem{prop}{Proposition}
\newtheorem{coroll}{Corollary}

\author{Jung Kyu {CANCI}}
\title{Finite orbits for rational functions
}
\date{}
\newcommand{\nc}{\newcommand}
\nc{\pro}{\mathbb{P}_1}
\nc{\po}{\mathbb{P}_1(K)}
\nc{\pq}{\mathbb{P}_1(\mathbb{Q})}
\nc{\ciclon}{(P_0, P_1,\ldots,P_{n-1})}
\nc{\orbit}{(P_{-m},\ldots,P_{-1},P_0,\ldots,P_{n-1})}
\nc{\cicloP}{\{P_i\}}
\nc{\tildep}{$\widetilde{a}^{\hspace{-0.65mm}\phantom{\tiny.}^{\tiny p}}$}
\nc{\spazio}{\hspace{0.7mm}}
\nc{\bdm}{\begin{displaymath}}
\nc{\edm}{\end{displaymath}}
\nc{\beq}{\begin{equation}}
\nc{\eeq}{\end{equation}}
\nc{\ind}{\hspace{-1.3mm}{\phantom{\_}}_}
\nc{\got}{\sffamily{p}\normalfont}
\nc{\ri}{R_S}
\nc{\ru}{R_S^\ast}
\nc{\qc}{\emph{quasi coprime}}
\nc{\p}{\mathfrak{p}}
\nc{\vvs}{\textquotedblleft}
\nc{\vvd}{\textquotedblright}
\nc{\dip}{\delta_{\p}}
\nc{\vap}{v_{\p}}
\begin{document}
\maketitle
\begin{abstract}\noindent Let $K$ be a number field and $\phi\in K(z)$ a rational function. Let $S$ be the set of all archimedean places of $K$ and all non-archimedean places associated to the prime ideals of bad reduction for $\phi$. We prove an upper bound for the length of finite orbits for $\phi$ in $\po$ depending only on the cardinality of $S$.
\end{abstract}
\section*{Introduction}
Let $K$ be a number field and $R$ its ring of integers. With every rational function $\phi\in K(z)$ we associate in the canonical way a rational map  $\Phi\colon\pro\to\pro$ defined over $K$. For every point $P\in\po$ we call its \emph{forward orbit under $\Phi$} (or simply \emph{orbit}) the set $O_{\Phi}(P)=\{\Phi^n(P)\mid n\in\mathbb{Z}_{\geq 0}\}$, where $\Phi^n$ is the $n$-th iterate of $\Phi$ and $\Phi^0(P)=P$. If $O_{\Phi}(P)$ is a finite set one says that $P$ is a \emph{pre-periodic} point for $\Phi$. This definition is due to the following fact: if $O_{\Phi}(P)$ is finite then there exist two integers $n\in\mathbb{Z}_{\geq 0}$ and $m\in\mathbb{Z}_{> 0}$ such that $\Phi^n(P)=\Phi^{n+m}(P)$. In this case one says that $\Phi^n(P)$ is a \emph{periodic} point for $\Phi$. If $m$ is the smallest positive integer with the above property, then one says that $m$ is the \emph{period} of $P$. If $P$ is a periodic point then its orbit is called a \emph{cycle}. 

It is not difficult to prove that every polynomial in $\mathbb{Z}[x]$ has cycles in $\mathbb{Z}$ of length at most $2$ and every finite orbit has cardinality at most $6$.
For a fixed finite set $S$ of valuations of $K$, containing all the archimedean ones, Narkiewicz in \cite{N.1} has shown that if $\Phi$ is a monic polynomial with coefficients in the ring of $S$-integers $R_S$ (see the definition at the beginning of the next section), then the length of its cycles in $K$ is bounded by a function $B(R_S)=C^{|S|(|S|+2)}$, for an absolute constant $C$. Note that the bound depends only on the cardinality of $S$. The value of $B(R_S)$ has been diminished by Pezda in \cite{P.1}. Indeed, the main result of Pezda \cite[Theorem 1]{P.1}, which concerns polynomial maps in local rings, combined with the estimate given in \cite[Theorem 4.7]{A.1} on the height of the $|S|$-th rational prime, gives rise to the following inequality
\beq\label{B_R} B(R_S)\leq \left[12|S|\log(5|S|)\right]^{2[K:\mathbb{Q}]+1}.\eeq
Later Narkiewicz and Pezda in \cite{N.P.1} extended \cite[Theorem 1]{P.1} to finite orbits so including pre-periodic points. By considering the limit in (\ref{B_R}) and the Evertse's bound proved in \cite{E.2} for the number of $S$-unit non-degenerate solutions to linear equations in three variables, the result of Narkiewicz and Pezda \cite[Theorem 1]{N.P.1} states that the length of a finite orbit in $K$ for a monic polynomial with coefficients in $R_S$ is at most
\bdm \frac{1}{3}\left[12|S|\log(5|S|)\right]^{2[K:\mathbb{Q}]+1}\left(31+2^{1031|S|}\right)-1.\edm

R. Benedetto has recently obtained a different bound, again for polynomial maps, but his bound also depends on the degree of the map. He proved in \cite{B.1} that if $\phi\in K[z]$ is a polynomial of degree $d\geq2$ which has bad reduction at $s$ primes of $K$, then the number of pre-periodic points of $\phi$ is at most $O(s\log s)$. The big-$O$ constant is essentially $(d^2-2d+2)/\log d$ for large $s$. Benedetto's proof relies on a detailed analysis of $\mathfrak{p}$-adic Julia sets. 

In the present paper we will generalize to finite orbits for rational maps the result of Narkiewicz and Pezda \cite{N.P.1} obtained for polynomial maps. We will study the same semigroup of rational maps studied in \cite{C.1}, namely: we fix an arbitrary finite set $S$ of places of $K$ containing all archimedean ones and consider the rational maps with good reduction outside $S$. We recall the definition of good reduction for a rational map at a non zero prime ideal $\p$ (for the details see \cite{M.S.2} or \cite{C.1}): a rational map $\Phi\colon\pro\to\pro$, defined over $K$, has good reduction at a prime ideal $\p$ if there exists a rational map $\tilde{\Phi}\colon\pro\to\pro$, defined over $K(\p)$, such that $\deg\Phi=\deg\tilde{\Phi}$ and the following diagram
\begin{center}
$\begin{diagram}
\node{\ \mathbb{P}_{1,K}}\arrow{e,t}{\Phi} \arrow{s,l}{\widetilde{\phantom{\Phi}}}
\node{\ \mathbb{P}_{1,K}}\arrow{s,r}{\widetilde{\phantom{\Phi}}}\\
\node{\ \mathbb{P}_{1,K(\mathfrak{p})}}\arrow{e,t}{\widetilde{\Phi}}\node{\ \mathbb{P}_{1,K(\mathfrak{p})}}
\end{diagram}$\end{center} is commutative, where $\ \tilde{}\ $ is the reduction modulo $\p$. In other words, an endomorphism $\Phi$ of $\pro$ defined over $K$ has good reduction at $\p$ if $\Phi$ can be written as $\Phi([X:Y])=[F(X,Y),G(X;Y)]$, where $F$ and $G$ are homogeneous polynomials of the same degree, with coefficients in the local ring $R_{\p}$ of $R$ at $\p$, and such that the resultant Res$(F,G)$ of polynomials $F$ and $G$ is a $\p$-unit in $R_{\p}$. Note that, from this definition, a rational map on $\po$ associated to a polynomial in $K[z]$ has good reduction outside $S$ if and only if its coefficients are $S$-integers and its leading coefficient is an $S$-unit.

In this paper we prove: 
\begin{teo}\label{Tp1}Let $K$ be a number field. Let $S$ be a finite set of cardinality $s$ of places of $K$ containing all the archimedean ones. There exists a number $c(s)$, depending only on $s$, such that the length of every finite orbit in $\po$, for rational maps with good reduction outside $S$, is bounded by $c(s)$.
We can choose $c(s)$ equal to 
\begin{align}\label{bound}\left[e^{10^{12}}(s+1)^8(\log(5(s+1)))^8\right]^{s}.\end{align} \end{teo}
The proof of Theorem \ref{Tp1} uses two non-elementary facts: the first is \cite[Corollary B]{M.S.1} where Morton and Silverman proved that if $\Phi$ is a rational map of degree $\geq 2$ which has bad reduction only at $t$ prime ideals of $K$ and $P\in\po$ is a periodic point with minimal period $n$, then the inequality
\beq\label{n} n\leq \left[12(t+2)\log(5(t+2))\right]^{4\left[K:\mathbb{Q}\right]}\eeq holds. The second one is the theorem proved by Evertse, Schlickewei and Schmidt in \cite{E.S.S.1} on the number of non-degenerate solutions $(u_1,\ldots,u_n)\in\Gamma$ to equation $a_1u_1+\ldots+a_nu_n=1$ where $\Gamma$ is a given subgroup of $(\mathbb{C}^\ast)^n$ of finite rank and the $a_i's$ are given non-zero elements of $K$. For $n=2$ and $a_1=a_2=1$ we use the upper bound proved by Beukers and Schlickewei in \cite{B.S.1}. The main point to obtain the estimate of Theorem \ref{Tp1} is the fact that the upper bounds in the theorems in \cite{E.S.S.1} and \cite{B.S.1} only depend on the rank of $\Gamma$. From Theorem \ref{Tp1} we easily deduce the following result concerning finite orbits for rational maps contained in a given finitely generated semigroup of endomorphisms of $\pro$:
\begin{coroll}\label{fgg}Let $\mathcal{F}$ be a finitely generated semigroup of endomorphisms 
of $\pro$ defined over a number field $K$. There exists a uniform upper bound $C$ which bounds the length of every finite orbit in $\po$ for any rational map in $\mathcal{F}$. Furthermore it is possible to give an explicit bound for $C$ in terms of a set of generators of $\mathcal{F}$. 
\end{coroll}

\bigskip

\emph{Acknowledgements.} The present work was written during the preparation of my Ph.D. thesis supervised by professor P. Corvaja. I would like to thank him for his useful suggestions. Also, I would like to thank prof. J.H. Silverman and prof. U. Zannier for their corrections and advice on my thesis. I am grateful to the referee for useful suggestions which helped to improve Theorem \ref{Tp1} and, in general, the presentation of this paper.
\section{Proofs}
In all the present paper we will use the following notation:
\begin{itemize}
\item[$K$]a number field;
\item[$R$]the ring of integers of $K$;
\item[$\mathfrak{p}$]a non zero prime ideal of $R$;
\item[$v_{\mathfrak{p}}$]the $\mathfrak{p}$-adic valuation on $R$ corresponding to the prime ideal $\mathfrak{p}$ (we always assume $v_\mathfrak{p}$ to be normalized so that $v_{\mathfrak{p}}(K^*)=\mathbb{Z}$); 
\item[$S$]a fixed finite set of places of $K$ of cardinality $s$ including all archimedean places;

\item[$R_S$] $\coloneqq\{x\in K \mid v_{\p}(x)\geq0 \ \text{for every prime ideal }\ \p\notin S\}$ the ring of $S$-integers;
\item[$R_S^\ast$] $\coloneqq\{x\in K^\ast\mid v_{\p}(x)=0 \ \text{for every prime ideal }\ \p\notin S\}$ the group of $S$-units.
\end{itemize}  
Let $P_1=\left[x_1:y_1\right]$ and $P_2=\left[x_2:y_2\right]$ be points in $\po$. Using the notation of  \cite{M.S.2} we will denote by  \begin{equation}\label{d_p}\delta_{\p}\,(P_1,P_2)=v_{\p}\,(x_1y_2-x_2y_1)-\min\{v_{\p}(x_1),v_{\p}(y_1)\}-\min\{v_{\p}(x_2),v_{\p}(y_2)\}\end{equation}the $\p$-adic logarithmic distance between the points $P_1,P_2$; note that  $\delta_{\p}\,(P_1,P_2)$ is independent of the choice of the homogeneous coordinates, i.e. it is well defined.
We will use the two following propositions contained in \cite{M.S.2}:
\begin{prop}\label{5.1}\emph{\cite[Proposition 5.1]{M.S.2}} 
\bdm \dip(P_1,P_3)\geq \min\{\dip(P_1,P_2),\dip(P_2,P_3)\}\edm
for all $P_1,P_2,P_3\in\po$.\vspace{-7.1mm}\begin{flushright}$\square$ \end{flushright}\end{prop}
\begin{prop}\label{5.2}\emph{\cite[Proposition 5.2]{M.S.2}} Let $\Phi\colon\mathbb{P}_1(K)\to\po$ be a rational map defined over $K$. Then
\bdm \dip (\Phi(P),\Phi(Q))\geq\dip(P,Q)\edm
for all $P,Q\in\po$ and all prime ideals $\p$ of good reduction for $\Phi$.\vspace{-7.1mm}\begin{flushright}$\square$ \end{flushright}
\end{prop}
With $(Q_{-m},\ldots,Q_0,\ldots,Q_{n-1})$ we always represent a finite orbit for a rational map $\Psi$ in which the $0$-th term $Q_0$ is a $n$-th periodic point for $\Psi$. Moreover, for all indexes $i\geq -m$, $Q_{i+1}=\Psi (Q_i)$ holds, bearing in mind that $Q_n=Q_0$.
We will use the following remark which is a direct consequence of the previous two propositions.
\begin{oss}\label{1+2} Let $(Q_{-m},\ldots,Q_0,\ldots,Q_{n-1})$ be a finite orbit in $\po$ for a rational map $\Psi$ with good reduction outside $S$; then for all integers $a,b$ with $-m\leq a\leq n-1,b\geq 0,k\geq0$ and for every prime ideal $\p\notin S$
\begin{multline*} \dip(Q_a,Q_{a+kb})\geq\\ \min\{\dip(Q_a,Q_{a+b}),\dip(Q_{a+b},Q_{a+2b}),\ldots,\dip(Q_{a+(k-1)b},Q_{a+kb})\}= \dip(Q_a,Q_{a+b}).\end{multline*}\end{oss}
\begin{proof} It is a direct application of the triangle inequality (Proposition \ref{5.1}) and Proposition \ref{5.2}. In fact the $b$-th iterate of $\Psi$ has good reduction at every prime ideal $\p\notin S$, therefore
\bdm \dip(Q_{a+lb},Q_{a+(l+1)b})=\dip(\Psi^b(Q_{a+(l-1)b}),\Psi^b(Q_{a+lb}))\geq\dip(Q_{a+(l-1)b},Q_{a+lb})\edm
for all indexes $0<l\leq k$.
\end{proof}

In the first version of this paper, in Theorem \ref{Tp1}, we proved an upper bound of the form $c(s,h)$ also depending on the class number $h$ of the ring $R_S$. Indeed we worked with a set $\mathbb{S}$ of places of $K$ containing $S$ such that the ring $R_{\mathbb{S}}$ was a principal ideal domain. From a simple inductive argument it results that it is possible to choose $\mathbb{S}$ such that $|\mathbb{S}|\leq s+h-1$. From some suitable applications of Proposition \ref{5.1} and Proposition \ref{5.2}, we obtained some equations in two and three $\mathbb{S}$-units and by using the upper bounds proved by Evertse in \cite{E.1} and \cite{E.2} we deduced a bound in Theorem \ref{Tp1}. Following the useful suggestions made by the anonymous referee we shall use, instead of the classical $S$-unit equation theorem, the refined result of Evertse, Schlickewei and Schmidt \cite{E.S.S.1} (and of Beukers and Schlickewei \cite{B.S.1} for $n=2$) leading to an upper bound in Theorem \ref{Tp1} depending only on the cardinality of $S$, even if $R_S$ is not a principal ideal domain. Now we state the last two quoted theorems and then we present the referee's suggestion to use these results.

Let $L$ be a number field. Let $(L^\ast)^n$ be the $n$-fold direct product of $L^\ast$, with coordinatewise multiplication $(x_1,\ldots,x_n)(y_1,\ldots,y_n)=(x_1y_1,\ldots,x_ny_n)$ and exponentiation $(x_1,\ldots,x_n)^t=(x_1^t,\ldots,x_n^t)$. We say that a subgroup $\Gamma$ of $(L^\ast)^n$ has rank $r$ if $\Gamma$ has a free subgroup $\Gamma_0$ of rank $r$ such that for every ${\bf x}\in \Gamma$ there is $m\in \mathbb{Z}_{>0}$ with ${\bf x}^m\in \Gamma_0$.\vspace{3 mm}

{\bf Theorem A} \cite{B.S.1} \emph{Let $L$ be a number field and let $\Gamma$ be a subgroup of $(L^\ast)^2$ of rank $r$. Then the equation 
\bdm x+y=1 \ \ \ \text{in $(x,y)\in\Gamma$}\edm
has at most $2^{8(r+1)}$ solutions.}\vspace{-7.1mm}\begin{flushright}$\square$ \end{flushright}\vspace{3 mm}

{\bf Theorem B} \cite{E.S.S.1} \emph{Let $L$ be a number field, let $n\geq 3$ and let $a_1,\ldots,a_n$ be non zero elements of $L$. Further, let $\Gamma$ be a subgroup of $(L^\ast)^n$ of rank $r$. Then the equation 
\bdm a_1x_1+\ldots+a_nx_n=1 \ \ \ \text{in $(x_1,\ldots,x_n)\in\Gamma$}\edm
has at most $e^{(6n)^{3n}(r+1)}$ solutions such that
$\sum_{\substack{i\in I}}a_ix_i\neq 0$ for each non empty subset $I\subset\{1,\ldots,n\}$}.\vspace{-8mm}\begin{flushright}$\square$ \end{flushright}\vspace{3 mm}

Let ${\bf a_1},\ldots,{\bf a_h}$ be a full system of representatives for the ideal classes of $R_S$. For each $i\in\{1,\ldots,h\}$ there is an $S$-integer $\alpha_i\in R_S$ such that 
\beq\label{alpha}{\bf a}_i^h=\alpha_iR_S.\eeq Let $L$ be the extension of $K$ given by 
\beq\label{lek} L=K(\zeta,\sqrt[h]{\alpha_1},\ldots, \sqrt[h]{\alpha_h})\eeq
where $\zeta$ is a primitive $h$-th root of unity. Of course if $h=1$ then $L=K$. Let us define the following subgroups of $L^\ast$
\bdm \sqrt{K^\ast}\coloneqq\{a\in L^\ast\mid \text{ $\exists $ $m\in\mathbb{Z}_{>0}$ with $a^m\in K^\ast$}\}\edm and
\bdm \sqrt{R_S^\ast}\coloneqq\{a\in L^\ast\mid \text{ $\exists $ $m\in\mathbb{Z}_{>0}$ with $a^m\in R_S^\ast$}\}.\edm
Let {\bf S} denote the set of places of $L$ lying above the places in $S$ and denote by $R_{\bf S}$ and $R_{\bf S}^\ast$ the ring of ${\bf S}$-integers and the group of ${\bf S}$-units, respectively, in $L$. By definition it is clear that $R_{\bf S}^\ast\cap \sqrt{K^\ast}= \sqrt{R_S^\ast}$ and so it follows that  $\sqrt{R_S^\ast}$ is a subgroup of $L^\ast$ of rank $s-1$. With the just stated notation we prove the following:

\begin{prop}\label{np}Let $L$ and ${\bf S}$ be as above. Let $\Phi$ be a rational map from $\pro$ to $\pro$ defined over $L$, having good reduction at all prime ideals outside ${\bf S}$. Let 
\beq\label{cn} \{P_{-m},\ldots,P_{-1},P_0\}\eeq 
be a set of $m+1$ distinct points of $\pro(L)$ such that $\Phi(P_i)=P_{i+1}$ for all $i\in\{-m,\ldots,-1\}$ and $\Phi(P_0)=P_{0}$. Further, suppose that $P_i=[x_i:y_i]$ for all indexes $i\in\{-m,\ldots,0\}$, where $x_i,y_i\in L$ such that
\begin{description}
\item{(1)} $x_0=0,y_0=1$;
\item{(2)} $x_iR_{\bf S}+y_iR_{\bf S}=R_{\bf S}$ for all indexes $i\in\{-m,\ldots,0\}$;
\item{(3)} $x_iy_j-x_jy_i\in\sqrt{K^\ast}$ for any distinct indexes $i,j\in\{-m,\ldots,0\}$.
\end{description}

Then $m< e^{10^{12}s}-2$.\end{prop}

The proof of this proposition will be a direct consequence of the following three lemmas.

\begin{lemma}\label{preliminare}With the same hypothesis of Proposition \ref{np}, let $P_{l-k},\ldots,P_{l-1},P_{l}$ be distinct points of the orbit (\ref{cn}) such that for every prime ideal $\p\notin {\bf S}$ 
\beq\label{l}\dip(P_{l-i},P_{0})=\dip(P_{l},P_{0})\ \ \text{for every index $0\leq i\leq k$}.\eeq Then $k< 2^{16^s}$.\end{lemma} 
\begin{proof} 
For every prime ideal $\p\notin {\bf S}$ and for any two indexes $k\geq i> j\geq0$ from Proposition \ref{5.1} and condition (\ref{l}) it follows that
\beq\label{d10} \dip(P_{l-i},P_{l-j})\geq \min\{\dip(P_{l-i},P_{0}),\dip(P_{l-j},P_{0})\}=\dip(P_{l},P_{0}).\eeq
Moreover, since $P_n=P_0$ for all $n\geq0$, by applying Remark \ref{1+2} to the orbit $(P_{-m},\ldots,P_{-1},P_0)$ with $a=l-i$, $b=i-j$ and $k=(m+1)$, where $m$ is the maximum integer such that $l-i+m(i-j)<0$, it follows that
\begin{align*} \dip(P_{l-i},P_{0})&\geq \min\{\dip(P_{l-i},P_{l-j}),\dip(P_{l-j},P_{l+i-2j}),\ldots ,\dip(P_{l-i+m(i-j)},P_0)\}\\&=\dip(P_{l-i},P_{l-j}).\end{align*}
By the last inequality, (\ref{l}) and (\ref{d10}) we have that
\beq\label{l-i}\dip(P_{l-i},P_{l-j})=\dip(P_{l},P_{0}).\eeq
Note that by condition \emph{(2)}-Proposition \ref{np} 
\beq\label{dv} \delta_{\p}(P_i,P_j)=v_{\p}(x_iy_j-x_jy_i)\eeq
for all indexes $i,j\in\{-m,\ldots,0\}$ and every prime ideal $\p\notin {\bf S}$.
Since $P_0=[0:1]$, from \emph{(3)}-Proposition \ref{np} it follows that $x_{l-i}\in\sqrt{K^\ast}$ and so, by (\ref{dv}) , condition (\ref{l}) is equivalent to $x_{l-i}R_{\bf S}=x_{l}R_{\bf S}$, for every index $0\leq i\leq k$.
Hence 
\beq\label{u_i}u_{l-i}\coloneqq\frac{x_{l-i}}{x_{l}}\in R_{\bf S}^\ast\cap \sqrt{K^\ast}= \sqrt{R_S^\ast}\eeq 
and $P_{l-i}=[x_{l}:y_{l-i}/u_{l-i}]$.
Furthermore, again from \emph{(3)}-Proposition \ref{np} combined with (\ref{l-i}) and (\ref{u_i}) we deduce that
\beq\label{sisu} u_{l-i,l-j}\coloneqq\frac{x_{l-i}y_{l-j}-x_{l-j}y_{l-i}}{x_lu_{l-i}u_{l-j}}=\frac{y_{l-j}}{u_{l-j}}-\frac{y_{l-i}}{u_{l-i}}\in R_{\bf S}^\ast\cap \sqrt{K^\ast}= \sqrt{R_S^\ast}\eeq
for all distinct indexes $i,j\in\{0,\ldots,k\}$.
In particular, either $k\in\{0,1\}$ or we have a system of three equations 
\beq\label{y_{l-i}}\left \{ \begin{array}{l}
y_l-y_{l-1}/u_{l-1}=u_{l-1,l}\\
y_l-y_{l-i}/u_{l-i}=u_{l-i,l}\\
y_{l-1}/u_{l-1}-y_{l-i}/u_{l-i}=u_{l-i,l-1}\end{array}\right..\eeq
The first one is obtained from (\ref{sisu}) substituting $j=0$ and $i=1$ and the two other ones with $j=0,j=1$ and $i$ an arbitrary index $k\geq i\geq 2$ (recall that $u_l=1$).\\
We deduce from (\ref{y_{l-i}}) the following linear relation:
\bdm u_{l-1,l}+u_{l-i,l-1}=u_{l-i,l},\edm
so $(u_{l-1,l}/u_{l-i,l},u_{l-i,l-1}/u_{l-i,l})\in \sqrt{R_S^\ast}\times \sqrt{R_S^\ast}$ is a solution of the equation $u+v=1$. Note that the group $\sqrt{R_S^\ast}\times \sqrt{R_S^\ast}$ has rank $2(s-1)$ therefore, by Theorem A (Beukers and Schlickewei \cite{B.S.1}) with $\Gamma=\sqrt{R_S^\ast}\times \sqrt{R_S^\ast}$, there are at most $2^{8(2s-2+1)}=2^{16s-8}$ possibilities for $(u_{l-1,l}/u_{l-i,l},u_{l-i,l-1}/u_{l-i,l})$. Now from (\ref{y_{l-i}}) it follows that
\bdm \frac{y_{l-i}}{u_{l-i}}=y_l-\frac{u_{l-i,l}}{u_{l-1,l}}u_{l-1,l}.\edm 
Thus the set of points $\left\{P_{l-i}=[x_{l}:y_{l-i}/u_{l-i}]\mid k\geq i\geq 2\right\}$ has cardinality bounded by $2^{16s-8}$ so $k\leq 2^{16s-8}+1<2^{16s}$.
\end{proof}
The next step is to prove an upper bound, which depends only on $s$, for the number of points $P_{-i}$ of (\ref{cn}) such that $x_{-i}R_{{\bf S}}\neq x_{-i+1}R_{{\bf S}}$. We need two lemmas. 

We say that a ${\bf S}$-integer $T$ is representable in \emph{two essentially different ways} as sum of two elements of $\sqrt{R_S^\ast}$ if there exist  
\beq\label{T} \text{$u_1,u_2,v_1,v_2\in \sqrt{R_S^\ast}$ such that $\{u_1,u_2\}\neq\{v_1,v_2\}$ and $T=u_1+u_2=v_1+v_2$}.\eeq

\begin{lemma}The cardinality of the set of non zero principal ideals of $R_{{\bf S}}$
\bdm \left\{T\cdot R_{{\bf S}}\mid \text{$T$ satisfies (\ref{T}})\right\}\edm is bounded by $e^{(18)^9(3s-2)}$.\end{lemma}
\begin{proof} Let $T\in R_{{\bf S}}/\{0\}$ be written as $T=u_1+u_2=v_1+v_2$ which satisfies the condition in (\ref{T}). Therefore the left term of equation
\bdm \frac{u_1}{v_1}+\frac{u_2}{v_1}-\frac{v_2}{v_1}=1\edm
has no vanishing subsums. Now, applying Theorem B (Evertse, Schlickewei and Schmidt \cite{E.S.S.1}) with $n=3$ and  $\Gamma=\sqrt{R_S^\ast}\times \sqrt{R_S^\ast}\times \sqrt{R_S^\ast}$ we obtain that the principal ideal 
\bdm T\cdot R_{{\bf S}}=v_1\left(1+\frac{v_2}{v_1}\right)\cdot R_{{\bf S}}\edm
has at most $e^{(18)^9(3s-2)}$ possibilities.\end{proof}
\begin{oss}\label{r1}\emph{By previous lemma, we can choose a set $\mathfrak{T}$ of ${\bf S}$-integers, with cardinality at most $e^{(18)^9(3s-2)}$, such that every non zero ${\bf S}$-integer with the property (\ref{T}) is representable as $uT$, where $u\in R_{{\bf S}}^*$ and $T\in\mathfrak{T}$}.\end{oss}
\begin{lemma}With the same hypothesis of Proposition \ref{np}, if there exist five distinct points  $P_{n_5}=[x_{n_5}:y_{n_5}],P_{n_4}=[x_{n_4}:y_{n_4}],P_{n_3}=[x_{n_3}:y_{n_3}],P_{n_2}=[x_{n_2}:y_{n_2}],P_{n_1}=[x_{n_1}:y_{n_1}]$ of the orbit (\ref{cn}), with ${n_5}<{n_4}<{n_3}<{n_2}<{n_1}<0$, then $x_{n_1}/x_{n_2}$ is a non zero ${\bf S}$-integer that is representable, in two essentially different ways, as sum of two elements of  $\sqrt{R_S^\ast}$.\end{lemma}
\begin{proof}Since $\Phi(P_0)=P_0=[0:1]$, from Proposition \ref{5.2},  considering $\Phi^{n_i-n_j}$, $P=P_{n_j}$ and $Q=P_0$, it follows that $x_{n_j}|x_{n_i}$ in $R_{{\bf S}}$ for all couple of integers $j\geq i$. Therefore there exist four non zero ${\bf S}$-integers $T_1,T_2,T_3,T_4$ such that 
\bdm x_{n_i}=T_ix_{n_{i+1}}\ \ \text{for all $i\in\{1,2,3,4\}$}\edm
and so for every couple of distinct indexes $1\leq i<j\leq5$ 
\beq\label{nj} x_{n_{i}}=T_i\cdot\ldots\cdot T_{j-1}x_{n_{j}}.\eeq
By Remark \ref{1+2} we have that 
\bdm \dip(P_{n_j},P_0)\geq\min\{\dip(P_{n_j},P_{n_i}),\dip(P_{n_i},P_{2n_i-n_j}),\ldots ,\dip(P_{m\cdot n_i-(m-1)n_j},P_0)\}=\dip(P_{n_i},P_{n_j})\edm for a suitable integer $m$, so it follows that $(x_{n_j}y_{n_i}-x_{n_i}y_{n_j})|x_{n_j}$ and by identity (\ref{nj})
\beq\label{egv} y_{n_i}-T_i\cdot\ldots\cdot T_{j-1}y_{n_{j}}=\frac{x_{n_j}y_{n_i}-x_{n_i}y_{n_j}}{x_{n_j}}\in R_{\bf S}^\ast\cap \sqrt{K^\ast}= \sqrt{R_S^\ast}.\eeq
(Recall that, by conditions \emph{(1)} and \emph{(3)} in the hypothesis of Proposition \ref{np}, the ${\bf S}$-integers  $x_{n_j}y_{n_i}-x_{n_i}y_{n_j}$ and $x_{n_j}$ belong to $\sqrt{K^\ast}$.)
From (\ref{egv}) we obtain:
\begin{align}
\label{a} &y_{n_1}-T_1y_{n_2}=v_1,\\
\label{b} &y_{n_2}-T_2y_{n_3}=v_2,\\
\label{c} &y_{n_1}-T_1T_2y_{n_3}=v_3,\\
\label{d} &y_{n_3}-T_3y_{n_4}=v_4,\\
\label{e} &y_{n_2}-T_2T_3y_{n_4}=v_5,\\
\label{f} &y_{n_1}-T_1T_2T_3y_{n_4}=v_6,\\
\label{g} &y_{n_2}-T_2T_3T_4y_{n_5}=v_7,\\
\label{h} &y_{n_1}-T_1T_2T_3T_4y_{n_5}=v_8,\\
\label{i} &y_{n_3}-T_3T_4y_{n_5}=v_9,\\
\label{j} &y_{n_4}-T_4y_{n_5}=v_{10}
,\end{align} where $v_i\in \sqrt{R_S^\ast}$ for all indexes $1\leq i\leq 10$.\\
From (\ref{a}), (\ref{c}) and (\ref{b}) we obtain
\beq\label{1} T_1=\frac{v_3}{v_2}-\frac{v_1}{v_2}.\eeq
From (\ref{f}), (\ref{a}) and (\ref{e}) we obtain
\beq\label{3} T_1=\frac{v_6}{v_5}-\frac{v_1}{v_5}.\eeq
From (\ref{h}), (\ref{a}) and (\ref{g}) we obtain
\beq\label{4} T_1=\frac{v_8}{v_7}-\frac{v_1}{v_7}.\eeq
Now we finish to proving that among (\ref{1}), (\ref{3}), (\ref{4}) there exist at least two distinct representations of $T_1$ as sum of two elements of $\sqrt{R_S^\ast}$.
From (\ref{g}), (\ref{b}) and (\ref{i}) we obtain that $T_2=\frac{v_7-v_2}{v_9}$; therefore $v_7\neq v_2$ and so

\beq\label{u1}\left\{\frac{v_3}{v_2},-\frac{v_1}{v_2}\right\}=\left\{\frac{v_8}{v_7},-\frac{v_1}{v_7}\right\} \Rightarrow -\frac{v_1}{v_2}=\frac{v_8}{v_7}.\eeq
From (\ref{g}), (\ref{e}) and (\ref{j}) we obtain that $T_2T_3=\frac{v_7-v_5}{v_{10}}$; therefore $v_7\neq v_5$ and so
\beq\label{u2}\left\{\frac{v_6}{v_5},-\frac{v_1}{v_5}\right\}=\left\{\frac{v_8}{v_7},-\frac{v_1}{v_7}\right\} \Rightarrow -\frac{v_1}{v_5}=\frac{v_8}{v_7}.\eeq
From (\ref{u1}) and (\ref{u2}) it follows that
\bdm \left\{\frac{v_3}{v_2},-\frac{v_1}{v_2}\right\}=\left\{\frac{v_6}{v_5},-\frac{v_1}{v_5}\right\}=\left\{\frac{v_8}{v_7},-\frac{v_1}{v_7}\right\} \Rightarrow -\frac{v_1}{v_2}=-\frac{v_1}{v_5}.\edm But this is not possible since, from (\ref{e}), (\ref{b}) and (\ref{d}), $T_2=\frac{v_5-v_2}{v_4}\neq 0$ holds.

\end{proof}
\begin{proof}[Proof of Proposition \ref{np}]
The set $\{P_{i_r},\ldots,P_{i_1},P_{-1}\}$ of all points $P_{-j}$ of the orbit (\ref{cn}) such that $x_{-j}R_{{\bf S}}\neq x_{-j+1}R_{{\bf S}}$ has cardinality equal to $r+1\leq 4+e^{(18)^9(3s-2)}$. Indeed, if such five points do no exist we have finish; otherwise for every index $i_{r-2}<i_t\leq i_1$ we apply the previous lemma with $n_1=-1,n_2=i_t,n_3=i_{r-2},n_4=i_{r-1},n_5=i_r$ obtaining that $x_{-1}x_{i_t}^{-1}=uT$  where $T\in\mathfrak{T}$ (the set chosen in Remark \ref{r1}) and $u$ is a suitable ${\bf S}$-unit. Therefore
\bdm P_{i_t}=[x_{-1}/T:uy_{i_t}].\edm  In this way we have proved that $r$ is bounded by $3+|\mathfrak{T}|$. Now, by Lemma \ref{preliminare}, it is clear that it is possible to choose as upper bound for $m$ the number
\bdm \left(4+e^{(18)^9(3s-2)}\right)\left(2^{16s}+1\right)< e^{(10)^{12}s}-2.\edm
\end{proof}

\begin{proof}[Proof of Theorem \ref{Tp1}] The bound (\ref{bound}) holds for finite orbit length for all rational maps of degree 1, i.e. automorphisms of $\mathbb{P}_1(K)$. Indeed every pre-periodic point for a bijection is a periodic point. Thus we have to study only the cycle lengths. If a point of $\mathbb{P}_1(K)$ is a periodic point for an automorphism $\Psi\in {\rm PGL}_2(K)$, with period $n\geq3$, then $\Psi^n$ is the identity map of $\mathbb{P}_1(K)$. The order of an element of ${\rm PGL}_2(K)$ is bounded by $2+4[K:\mathbb{Q}]^2$, so from $2s\geq [K:\mathbb{Q}]$ it results that 
\beq\label{aut}n\leq 2+16s^2<c(s).\eeq

Now we consider rational maps of degree $\geq 2$ with good reduction outside $S$. We reduce to the hypothesis of Proposition \ref{np}. Let $(Q_{-l},\ldots,Q_{-1},Q_0,\ldots,Q_{n-1})$ be a finite orbit in $\po$, for a rational map $\Psi\colon\pro\to\pro$ defined over $K$ with good reduction outside $S$, including $(Q_0,\ldots,Q_{n-1})$ as a cycle for $\Psi$. We can associate a finite orbit in which the cycle consists of one single point (i.e. a fixed point). Indeed, the tuple $(Q_{-\left\lfloor \frac{l}{n}\right\rfloor n},\ldots,Q_{-n},Q_O)$ is an orbit for $\Psi^n$ and $Q_0$ is a fixed point. We set $m\coloneqq\left\lfloor \frac{l}{n}\right\rfloor$. Of course $\Psi^n$ can be viewed as an endomorphism of $\pro$ defined over $L$ (the extension of $K$ defined in (\ref{lek})). For every index $i\in\{0,\ldots,m\}$, let $Q_{-i\cdot n}=[t_i:l_i]$ be a representation of $Q_{-i\cdot n}$ in $S$-integral homogeneous coordinates. Recall that $\{{\bf a_1},\ldots,{\bf a_h}\}$ is a full system of representatives for the ideal classes of $R_S$ and that the $\alpha_i's$ are the $S$-integers verifying (\ref{alpha}). Let ${\bf b}_i\in\{{\bf a}_1,\ldots,{\bf a}_h\}$ be the representative of the ideal $t_iR_S+l_iR_S$. Let $\beta_i\in\{\alpha_1,\ldots,\alpha_h\}$ be such that ${\bf b}_i^h=\beta_iR_S$. Hence there exists $\lambda_i\in K^\ast$ satisfying $(t_iR_S+l_iR_S)^h=\lambda_i^h\beta_iR_S$. As suggested by the referee, in $L$, we define 
\beq\label{skc} t_i^\prime\coloneqq \frac{t_i}{\lambda_i\sqrt[h]{\beta_i}},\ \ \ l_i^\prime\coloneqq \frac{l_i}{\lambda_i\sqrt[h]{\beta_i}}.\eeq
It is clear that $t_i^\prime,l_i^\prime$ are elements of $\sqrt{K^\ast}$ such that 
\beq\label{qpid}(t_i^\prime R_{\bf S}+l_i^\prime R_{\bf S})=R_{\bf S}.\eeq
Furthermore, for any two distinct indices $i,j$
\bdm (t_i^\prime l_j^\prime-t_j^\prime l_i^\prime)^h=\frac{(t_il_j-t_jl_i)^h}{\lambda^h_i\lambda_j^h\beta_i\beta_j}\in K^\ast.\edm

By (\ref{qpid}) with $i=0$ there exist $r_0,s_0\in R_{\bf S}$ such that $r_0t^{\prime}_0+s_0l^{\prime}_0=1$. Define the matrix 
\bdm A=\begin{pmatrix}
	l^{\prime}_0 & -t^{\prime}_0\\
	r_0 & s_0
\end{pmatrix}\edm
and further define $x_i,y_i$ by 
\bdm  \begin{pmatrix}
	x_i\\
	y_i
\end{pmatrix}=A\begin{pmatrix}
	t_i^\prime\\
	l_i^\prime
\end{pmatrix}\edm
for all $i\in\{0,\ldots,m\}$. If we now set $P_i\coloneqq[x_i:y_i]$ for all $i\in\{0,\ldots,m\}$ and $\Phi\coloneqq[A]\circ\Psi^n\circ[A]^{-1}$, where $[A]$ is the automorphism of $\pro$ induced by $A$, then, by Proposition \ref{np}, $m=\left\lfloor \frac{l}{n}\right\rfloor<e^{(10)^{12}s}-2$ and so $l<n(e^{(10)^{12}s}-1)$. Therefore the orbit $(Q_{-l},\ldots,Q_{-1},Q_0,\ldots,Q_{n-1})$ has cardinality bounded by $ne^{(10)^{12}s}$. Since in $S$ there are at most $s-1$ prime ideals and $2s\geq[K:\mathbb{Q}]$, the inequality (\ref{bound}) becomes 
\bdm n\leq \left[12(s+1)\log(5(s+1))\right]^{8s}\edm 
and so the theorem is proved.

\end{proof}

\begin{proof}[Proof of Corollary \ref{fgg}]
Choose a finite set of generators of $\mathcal{F}$. Each of these generators has at most finitely many prime ideals of $R$ of bad reduction. So there is a finite set $S$ of prime ideals such that each of the chosen generators, and therefore each of the elements of $\mathcal{F}$, has good reduction outside $S$. We conclude by applying Theorem \ref{Tp1}.

\end{proof}

\bigskip

\noindent Jung Kyu CANCI\\
Dipartimento di Matematica e Informatica\\
Università degli Studi di Udine\\
via delle Scienze, 206\\
33100 Udine (ITALY)\\
E-mail: {\sf canci@dimi.uniud.it}

\end{document}